\documentclass{amsart}
\usepackage{amssymb}
\newtheorem{theorem}{Theorem}[section]

\newtheorem{lemma}[theorem]{Lemma}
\newtheorem{claim}[theorem]{Claim}
\newtheorem{fact}[theorem]{Fact}
\newtheorem{corrolary}[theorem]{Corrolary}
\newtheorem{problem}{Problem}[]
\theoremstyle{definition}
\newtheorem{definition}[theorem]{Definition}

\theoremstyle{remark}

\numberwithin{equation}{section}

\author[Yona Cherniavsky]{Yona Cherniavsky}
\address{Department of Mathematics, Bar Ilan University, Ramat Gan 52900, Israel}
\email{cherniy@math.biu.ac.il , chrnvsk@gmail.com}

\author{Mishael Sklarz}

\title[Conjugacy in Permutation Representations of $S_n$]{Conjugacy in Permutation Representations of the Symmetric Group}

\subjclass[2000]{Primary 20C30; Secondary 05E15}
\thanks{Supported in part by a grant from 
the Israel Science Foundation.}
\keywords{conjugacy classes, symmetric group, permutation
representations, characters, fixed points}

\newcommand{\bF}{\mathbb{F}}
\newcommand{\bC}{\mathbb{C}}
\newcommand{\set}[1]{\left\{#1\right\}}

\newcommand{\gen}[1]{\left\langle#1\right\rangle}
\DeclareMathOperator{\chr}{char}
\DeclareMathOperator{\tr}{trace}
\DeclareMathOperator{\fix}{fix}

\begin{document}

\begin{abstract}
Although the conjugacy classes of the general linear group are
known, it is not obvious (from the canonic form of matrices) that
two permutation matrices are similar if and only if they are
conjugate as permutations in the symmetric group, i.e. that
conjugacy classes of $S_n$ do not unite under the natural
representation. We prove this fact, and give its application to
the enumeration of fixed points under a natural action of
$S_n\times S_n$. We also consider the permutation representations
of $S_n$ which arise from the action of $S_n$ on ordered tuples
and on unordered subsets, and classify which of them unite
conjugacy classes and which do not.
\end{abstract}

\maketitle

\section{Introduction}
In this paper we study the action of $S_n$ on ordered $k$-tuples.
Denote by $\rho_k^{\bF}$ the corresponding permutation
representation over an arbitrary field $\bF$. The following
problem was presented to us by A. Lubotzky and Y. Roichman.

\begin{problem} For which $1\le k\le n$ does the following hold:

For any two permutations $\pi,\sigma\in S_n$, $\rho_k^{\bF}(\pi)$
is conjugate to $\rho_k^{\bF}(\sigma)$ in $GL(n,{\bF})$ if and
only if $\pi$ and $\sigma$ are conjugate in $S_n$.
\end{problem}

This problem arises in the enumeration of invertible matrices with
respect to a certain natural action of $S_n\times S_n$, see
\cite{yona} and Section \ref{alpha} below. Our Theorem~4.5 is used
by C. de la Mora and P. Wojciechowski~\cite{MW} in enumerating of
multiplicative bases of matrix algebras.

\medskip

For $k=n$, a negative solution to Problem 1 was essentially given
by Burnside \cite{burnside}. In Section \ref{sec-1} it is shown
that for $k=1$ the answer is positive. A full solution is given in
Section \ref{sec-2}:  We find that $\rho_1^\bF$ and $\rho_2^\bF$
do not unite any classes, that $\rho_3^\bF$ unites classes only
when $n$ is even, and that $\rho_k^\bF$ for $k \geq 4$ always
unites some classes. In Subsection~\ref{subsets} we try to
generalize these results to the permutation representations of
$S_n$ arising from the action of $S_n$ on (unordered) $k$-subsets
of $\{1,2,...,n\}$: it appears that the representation arising
from the action of $S_n$ on {\it all} subsets of $\{1,2,...,n\}$
does not unite any classes (Theorem~\ref{allsubsets});
 when $n$ is odd, the action of $S_n$ on the set
of even-sized subsets of $\{1,2,...,n\}$ does not unite classes
(Theorem~\ref{evensizedsubsets}) and finally, the action of $S_n$
on the set of odd-sized subsets of $\{1,2,...,n\}$ does not unite
classes and this does not depend on n's parity
(Theorem~\ref{oddsizedsubsets}). These results do not depend on
the choice of the field $\bF$. Finally, our results are applied in
Section \ref{alpha} to the enumeration of fixed points of a
natural action of $S_n\times S_n$ on invertible matrices.

\section{The Natural Representation of $S_n$}\label{sec-1}
There is an obvious embedding of $S_n$ in $GL(n,\bF)$ where $\bF$
is any field.  Consider  a permutation $\pi\in S_n$ as an $n\times
n$ matrix obtained from the identity matrix by permutations of the
rows. More explicitly: for every permutation $\pi \in S_n$ we
identify $\pi$ with the matrix:
\begin{center}
$$[\pi]_{i,j} =  \left\{ \begin{array}{cc}
{1}      &  { i=\pi(j)} \\
{0 }          & \text{otherwise}
\end{array}
\right. $$
\end{center}
This representation can also be realized as the permutation
representation which is obtained from the natural action of $S_n$
on $\set{1,2,\dots,n}$ defined by $\pi\cdot i = \pi(i)$.

Our first result is that this representation does not unite
conjugacy classes of $S_n$.  We shall use the following well known
fact:

\begin{fact}\label{split}
If $\sigma$ is a cycle of length $n$, then $\sigma^k$ consists of
$(n,k)$ cycles, each of length $n/(n,k)$.\footnote{We use $(n,k)$
to denote the greatest common divisor of $n$ and $k$. }
\end{fact}

\begin{theorem}\label{similar}
  The conjugacy classes of $S_n$ do not
unite in $GL(n,\bF)$. In other words, if $\pi$ and $\sigma$ are
permutations with similar matrices in $GL(n,\bF)$, then they are
conjugate in $S_n$ too.
\end{theorem}

\begin{proof} Let $\pi$ and $\sigma$ be permutations which
are similar as matrices. First of all, we note that for any $k$,
$\pi^k$ and $\sigma^k$ are also similar.

\noindent\textbf{Case 1:} $\chr(\bF)=0$ or at least $\chr(\bF)>n$.

Each cycle of length $k$ in $\pi$ contributes the term $x^k-1$
into the characteristic polynomial of the permutation matrix.
Under the above restriction on $\chr(\bF)$ it seems reasonable
that the cycle structure of a permutation can be recovered from
the characteristic polynomial of the corresponding permutation
matrix. However, our proof utilizes the trace of the permutation
matrix and the traces of its powers.

Denote by $c_d(\pi)$ the number cycles with length equal to $d$ in
$\pi$. We shall use induction on $d$ to prove that
$c_d(\pi)=c_d(\sigma)$, for all $d$, and this will show that $\pi$
and $\sigma$ are conjugate.

Since $\pi$ and $\sigma$ are similar as matrices, we have
$\tr(\pi)=\tr(\sigma)$.  However, the trace function counts the
$1$'s on the diagonal (here we use the restriction on
$\chr(\bF)$), and each such $1$ corresponds to a fixed point of
the permutation, so $\tr(\pi)=c_1(\pi)$. Therefore,
$c_1(\pi)=c_1(\sigma)$, i.e. $\pi$ and $\sigma$ have the same
number of fixed points.  This is the base of our induction.

Now let $d$ be an arbitrary number, and suppose that
$c_k(\pi)=c_k(\sigma)$ for all $k<d$. From Lemma \ref{split} it
follows that a $k$-cycle in $\pi$ ends up as a product of $k$
$1$-cycles in $\pi^d$ if and only if $k$ divides $d$. Therefore,
we can conclude that
$$\tr(\pi^d)=\sum_{k\mid d} k\cdot c_k(\pi) = d\cdot c_d(\pi) + \sum_{k\mid d, k<d} k\cdot c_k(\pi).$$
Now, by our induction hypothesis, for all proper divisors $k \mid
d$ we have $c_k(\pi)=c_k(\sigma)$.  On the other hand,
$\tr(\pi^d)=\tr(\sigma^d)$.   This implies that
$c_d(\pi)=c_d(\sigma)$, and completes the induction argument.

We have shown that $\pi$ and $\sigma$ have the same cycle
structure, so they are conjugate as permutations.

\noindent\textbf{Case 2:}  With no restriction on $\chr(\bF)$

In this case, the trace of a permutation matrix no longer gives
the number of fixed points of the permutation, so a more devious
route is necessary.

Note that in this case it is impossible to recover the cycle
structure of a permutation from the characteristic polynomial of
the corresponding permutation matrix: for example, if
$\chr(\bF)=2$ we have $x^4+1=(x^2+1)^2=(x+1)^4$, i.e. one cycle of
length 4, two cycles of length 2 and four cycles of length 1 have
the same characteristic polynomial.

\begin{claim}\label{cl1}
Denote by $m(\pi)$ the number of cycles in $\pi$. $m(\pi)$ is
invariant for similar matrices, i.e., $m(\pi)=m(\sigma)$
\end{claim}
\begin{proof}
The number $1$ is an eigenvalue of all permutation matrices,
because $(1,1,\ldots,1)^T$ is an eigenvector.
Since $\pi$ acts on vectors by permuting their coordinates,
solving the equation $\pi\underline{x}=\underline{x}$ (which is
the equation to find eigenvecors of the eigenvalue 1) we have the
following: each cycle $(i_1,i_2,\ldots,i_k)$ in $\pi$ gives us $k$
equations on corresponding entries of $\underline x$
$$
x_{i_k}=x_{i_1}\,\,,\,\,x_{i_1}=x_{i_2}\,\,,\,\,\ldots\,\,,\,\,x_{i_{k-1}}=x_{i_k}
$$
which means that all these entries are equal. Thus we have that
the number of free variables, which is the dimension of
1-eigenspace, equal to the number of cycles in $\pi$. Noting that
the dimension of an eigenspace is invariant for similar matrices
proves our point.
\end{proof}

\begin{claim}\label{cl2}
Denote by $m_d(\pi)$ the number of cycles in $\pi$ with length
divisible by $d$.   $m_d(\pi)$ is also invariant for similar
matrices, i.e., $m_d(\pi)=m_d(\sigma)$
\end{claim}
\begin{proof}
First, take $d=2$ as an example.  Any even-length cycle in $\pi$
splits into two cycles in $\pi^2$, while the odd cycles stay
unsplit.  This means that
$m_2(\pi)=m(\pi^2)-m(\pi)=m(\sigma^2)-m(\sigma)=m_2(\sigma)$. More
generally, for any prime number $p$, the cycles of length
divisible by $p$ split into $p$ separate cycles, while cycles of
length non-divisible by $p$ remain unsplit.  Therefore,
$$m_p(\pi)=\frac{m(\pi^p)-m(\pi)}{p-1}\,\,,$$ which implies that
$m_p(\pi)=m_p(\sigma)$ by Claim \ref{cl1}.

Now turn to the proof in general case.

We shall proceed by induction on the length of the prime
decomposition of $d$. First, for $d=1$, $m_d(\pi)=m(\pi)$, and our
claim follows from Claim \ref{cl1}.

Suppose that $d=pt$, with $p$ prime.  An appeal to Lemma
\ref{split} shows that cycles in  $\pi$ of length divisible by
$pt$ split into $p$ separate cycles in $\pi^p$, each of length
divisible by $t$.  The cycles in $\pi$ whose length is
non-divisible by $pt$ either split into cycles of length
non-divisible by $t$ or don't split at all.  Working backwards, a
cycle in $\pi^p$ of length divisible by $t$ either comes from a
cycle in $\pi$ that was divisible by $pt$ (in which case it will
be one of $p$ such cycles), or comes from a cycle in $\pi$ of
length divisible by $t$ but \emph{not} by $p$.

If $t$ is prime to $p$, the number of cycles in $\pi$ divisible by
$t$ but not by $p$ is $(m_t(\pi)-m_{pt}(\pi))$, so we have
$m_t(\pi^p)=p\cdot m_{pt}(\pi) + (m_t(\pi)-m_{pt}(\pi))$. On the
other hand, if $t$ is divisible by $p$, there are no cycles
divisible by $t$ but not by $p$, so we have $m_{pt}(\pi)=p\cdot
m_t(\pi^p)$. In summary, we have shown that
$$m_{pt}(\pi)=\left\{ \begin{array}{ll}  \frac{m_t(\pi^p)-m_t(\pi)}{p-1} & (p,t)=1  \\
                     \frac{m_t(\pi^p)}{p}            & (p,t)>1
              \end{array} \right.
$$
Now, $t$ has a shorter prime decomposition, so by our induction
hypothesis, $m_t$ is invariant, so $m_t(\pi) = m_t(\sigma)$, and
$m_t(\pi^p) = m_t(\sigma^p)$. This shows that
$m_n(\pi)=m_n(\sigma)$, as was to be shown.

A short example is in place.  Suppose $\pi$ has cycle
decomposition $[3^2,6^2,9,12^2]$.  Then $\pi^2$ has decomposition
$[3^6,6^4,9]$, and $\pi^3$ decomposes as $[1^6,2^6,3^3,4^6]$. On
the one hand, we have $$m_{12}(\pi) = m_{3\cdot 4}(\pi) =
(m_4(\pi^3)-m_4(\pi))/(3-1) = (6-2)/2 = 2.$$ On the other hand,
$$m_{12}(\pi) = m_{2\cdot 6}(\pi) = m_6(\pi^2)/2 = 4/2 = 2,$$ and
we see that both calculations give us the expected results.
\end{proof}

\begin{claim}\label{cs}
$c_d(\pi)$ (the number of cycles of exact length $d$) is also
invariant, so $c_d(\pi) = c_d(\sigma)$.
\end{claim}
\begin{proof}
Fix some permutation $\pi$, and let $I_d$ denote the set of cycles
of length divisible by $d$ in $\pi$.  Evidently $m_d(\pi)=|I_d|$.
Also, we have $I_d \cap I_k = I_{[d,k]}$, so $|I_d \cap I_k| =
m_{[d,k]}(\pi)$.  The same goes for intersections of more than two
such sets.  Now, it is also clear that $c_d(\pi)=|I_d - \cup_{k>1}
I_{dk}| = |I_d|-|\cup_{k>1} I_{dk}|$. Using the
inclusion-exclusion principal, we end up with a sum
$$c_d(\pi) = |I_d| - \sum_{k>1} |I_{dk}| + \sum_{k_1,k_2 > 1} |I_{dk_1}\cap I_{dk_2}| - \dots$$
and as we saw, all the summands are in fact of the form
$m_k(\pi)$, for some $k$. In the previous claim we showed that all
the $m_k$ are invariant, so $c_d$ is invariant too, and our claim
is proved.

Again, a short example might clarify things.  Using $\pi$ from the
example above, with cycle type $[3^2,6^2,9,12^2]$, we have
\begin{align*}
c_3(\pi)& = |I_3| - (|I_6|+|I_9|+|I_{12}|) + (|I_6\cap I_9| + |I_6\cap I_{12}| + |I_9 \cap I_{12}|) \\
                &{\phantom{=}} - |I_6\cap I_9\cap I_{12}| \\
        & = |I_3| - (|I_6|+|I_9|+|I_{12}|) + (|I_{18}| + |I_{12}| + |I_{36}|) - |I_{36}|  \\
        & = m_3 - (m_6 + m_9 + m_{12}) + (m_{18} + m_{12} + m_{36}) - m_{36} \\
        & = 7 - ( 4 + 1 + 2) + (0+2+0) - 0 = 2
\end{align*}
and that indeed is the correct result.
\end{proof}
Summing up, Claim \ref{cs} means that permutation matrices $\pi$
and $\sigma$ which are similar as matrices have the same cycle
structure, i.e. are conjugate in the symmetric group and thus the
proof of Theorem \ref{similar} is completed.
\end{proof}

\section{Other Permutation Representations}\label{sec-2}

As we mentioned above, the natural representation of $S_n$ can be
realized as the permutation representation obtained from the
natural action of $S_n$ on $\set{1,2,\dots, n}$, defined by
$\pi\cdot i = \pi(i)$.  We now consider the permutation
representation $\rho_k^\bF$, which arise from the action of $S_n$
on $k$-tuples, where the action is defined by $\pi\cdot(i_1,
\dots, i_k) = (\pi(i_1), \dots, \pi(i_k))$.  These representations
form an interpolation between the natural representation
($1$-tuples) and the regular representation ($n$-tuples).

\subsection{Representations Over the Complex Field}
In Section \ref{sec-1} we proved that the natural representation
of $S_n$ does not unite conjugacy classes.  On the other hand, it
is well known \cite{burnside} that the regular representation of
$S_n$ (indeed, of any group) unites all elements of equal order.
The natural representation can be seen as the permutation
representation obtained from the natural action of $S_n$ on the
set $\set{(1),(2),\dots, (n)}$ of $1$-tuples.  On the other hand,
the regular representation can be seen as the permutation
representation which arises from the action of $S_n$ on all $n!$
ordered $n$-tuples of numbers from $\set{1,2,\dots,n}$.  In this
section we wish to address the representations in between:  the
representation arising from the action of $S_n$ on pairs,
triplets, etc. and to see where the representations start uniting
conjugacy classes.

We begin with a general result, which holds true for any
representation of any finite group.

\begin{theorem}\label{general}
Let $G$ be a group, and $\sigma,\tau\in G$. Let $T:G \to
GL(d,\bC)$ a representation of $G$, with character $\chi$. Then
$T(\sigma)\sim T(\tau)$ as matrices if and only if
$\chi(\sigma^k)=\chi(\tau^k)$ for all $k$.
\end{theorem}
\begin{proof}
Obviously, if $T(\sigma)\sim T(\tau)$, then $T(\sigma^k)\sim
T(\tau^k)$ for any $k$, so $\chi(\sigma^k)=\tr(T(\sigma^k)) =\tr
(T(\tau^k))=\chi(\tau^k)$.

Now suppose that $\chi(\sigma^k)=\chi(\tau^k)$ for all $k$. Denote
$s = |\sigma|$ and $t=|\tau|$, and let $m=[s,t]$.  Now consider
$C=\gen{x}$, a cyclic group or order $m$, and define two
representations of $C$:  $T_\sigma(x)=T(\sigma)$ and
$T_\tau(x)=T(\tau)$.  It is easily seen that these are indeed well
defined representations. The hypothesis now reads that these two
representations have the same character, and so they must be
similar \cite{curtis}.  This implies that $T(\sigma) \sim
T(\tau)$.
\end{proof}

Note that the fact that the regular representation unites all
elements of equal order can be derived from this theorem:  If
$\chi$ is the character of the regular representation, then
$$\chi(\sigma^k)=
  \begin{cases}
    |G| & |\sigma| \mid k, \\
     0 &  \text{otherwise} \\
  \end{cases}$$
so obviously $\chi(\sigma^k)=\chi(\tau^k)$ for all $k$ if and only
if $\sigma$ and $\tau$ have the same order.

We also note that representations which are not faithful trivially
unite all the conjugacy classes in their kernel with $1$.  We may
therefore omit them from our discussion from now on. Faithful
representations can only unite two classes if their elements have
the same order.

The criterion which we just presented is still rather complicated
to use for general groups, but it can be simplified in our case,
because of the following simple fact.

\begin{fact}
Let $\sigma\in S_n$, with $|\sigma|=m$.
\begin{itemize}
\item If $k$ is relatively prime to $m$.  Then $\sigma^k \sim
\sigma$.

\item For any $k$, $\sigma^k \sim \sigma^{(m,k)}$.
\end{itemize}
\end{fact}
%

\begin{claim}\label{sn-general}
Let $T:S_n\to GL(d,\bC)$ be a representation of the symmetric
group, with character $\chi$, and $\sigma, \tau \in S_n$ elements
of order $m$. Then $T(\sigma)\sim T(\tau)$ if and only if
$\chi(\sigma^k)=\chi(\tau^k)$ for $k\mid m$.
\end{claim}
\begin{proof}
If $T(\sigma)\sim T(\tau)$, then obviously
$\chi(\sigma^k)=\chi(\tau^k)$ for any $k$, and in particular for
divisors of $m$.  On the other hand, suppose
$\chi(\sigma^k)=\sigma(\tau^k)$ for $k\mid m$.  Now, for any $k$,
$\sigma^k\sim \sigma^{(k,m)}$, and similarly for $\tau$.
Therefore, using the fact that $(k,m)\mid m$,
$\chi(\sigma^k)=\chi(\sigma^{(k,m)}) = \chi(\sigma^{(k,m)}) =
\chi(\tau^k)$.  This holds for all $k$, and therefore, by Theorem
\ref{general}, $T(\sigma)\sim T(\tau)$.
\end{proof}

\begin{corrolary}
If $\sigma$ and $\tau$ are of prime order $p$, then $T(\sigma)\sim
T(\tau)$ if and only if $\chi(\sigma)=\chi(\tau)$.
\end{corrolary}

\begin{definition}
Let $\sigma, \tau \in S_n$ be elements of equal order $m$, such
that if $k\neq 1$ and $k\mid m$ then $\sigma^k \sim \tau^k$.  Then
directly from \ref{sn-general} it follows that $T(\sigma)\sim
T(\tau)$ if and only if $\chi(\sigma)=\chi(\tau)$. We call such
elements \emph{almost similar}.  In fact, it is sufficient to
require that $\sigma^p\sim \tau^p$ for all prime divisors of $m$.
\end{definition}

We next show that almost similar elements are typical examples of
elements that are united by representations, in the following
sense:

\begin{theorem}
Let $T:S_n\to GL(d,\bC)$ be a representation.  If $T$ unites some
two conjugacy classes, then there must exist a pair of almost
similar elements which it unites.
\end{theorem}
\begin{proof}
Suppose that $\sigma$ and $\tau$ are non-similar elements such
that $T(\sigma)\sim T(\tau)$.  Obviously $T(\sigma^k)\sim
T(\tau^k)$, so it suffices to show that there exists some $k$ for
which $\sigma^k$ and $\tau^k$ are almost similar but non-similar.
We show this by induction on the number of primes in the prime
decomposition of $m$, the order of $\sigma$ and $\tau$.  If $m$ is
of length $1$, i.e. $m$ is prime, then  $\sigma$ and $\tau$
already are almost similar.  In general, if $\sigma^p\sim \tau^p$
for all prime divisors $p\mid m$, then $\sigma$ and $\tau$ are
already almost similar.  Otherwise, there exists some prime $p
\mid m$ such that $\sigma^p \not\sim \tau^p$. The order of
$\sigma^p$ and $\tau^p$ is $m/p$, which has shorter prime
decomposition than $m$, so by the induction hypothesis there
exists some $k$ for which $\sigma^{pk}$ and $\tau^{pk}$ are almost
similar.
\end{proof}

Having proved this, we now have a criterion to check whether a
representation unites classes:  It is sufficient to show that all
pairs of almost similar elements remain non united, i.e. that the
character of the representation takes different values on them.

\begin{lemma}\label{almost-sim}
Let $\sigma$ and $\tau$ be almost similar permutations in $S_n$,
of order $m$.
\begin{enumerate}
\item If $\fix(\sigma)=\fix(\tau)$, then $\sigma \sim \tau$.

\item For any prime $p\mid m$, $p\mid (\fix(\sigma)-\fix(\tau))$.
\end{enumerate}
\end{lemma}
\begin{proof} \hfill
\begin{enumerate}
\item Suppose that $\fix(\sigma)=\fix(\tau)$.  We shall show that
$c_k(\sigma)=c_k(\tau)$, for all $k\mid m$, by induction on the
length of $k$'s prime decomposition. We have
$c_1(\sigma^k)=c_1(\tau^k)$. If $k$ is prime, then
$c_1(\sigma^k)=kc_k(\sigma)+c_1(\sigma)$, so after clearing sides,
we have $c_k(\sigma)=c_k(\tau)$.  In the general case, we have
$c_1(\sigma^k)=kc_k(\sigma) + \sum_{t\mid k, t\neq k}
tc_t(\sigma)$.  By the induction hypothesis, the summands under
the summation sign are equal for $\sigma$ and $\tau$, and so is
the left hand side.  This implies that $c_k(\sigma)=c_k(\tau)$.

\item Let $p\mid m$.  Then $c_1(\sigma^p)=c_1(\tau^p)$, so
$pc_p(\sigma)+c_1(\sigma)=pc_p(\tau)+c_1(\tau)$.  Clearing sides
gives $c_1(\sigma)-c_1(\tau)=p(c_p(\tau)-c_p(\sigma))$.
\end{enumerate}
\end{proof}

We are now in a position to re-prove the theorem on the natural
representation from Section \ref{sec-1} over $\bC$:

\begin{theorem}
The natural representation $\rho_1^\bC: S_n\to GL(n, \bC)$ does
not unite conjugacy classes.
\end{theorem}
\begin{proof}
In this case $\chi(\sigma)=c_1(\sigma)$.  Suppose $\sigma$ and
$\tau$ are almost similar and $\chi(\sigma)=\chi(\tau)$.
 This means that $c_1(\sigma)=c_1(\tau)$, and Lemma \ref{almost-sim}
 implies that $\sigma \sim \tau$.
\end{proof}

Using a similar technique, we can go further:
\begin{theorem}
The representation $\rho_2^\bC$ of $S_n$ on ordered pairs does not
unite classes either.
\end{theorem}
\begin{proof}
In this case $\chi(\sigma) = c_1(\sigma)(c_1(\sigma)-1)$, because
each fixed point of $T(\sigma)$ comes from an ordered pair of
fixed points of $\sigma$.  Suppose there exists a pair $\sigma$
and $\tau$ of almost similar elements of order $m$ which are
united by $T$, i.e. $\chi(\sigma)=\chi(\tau)$. For convenience,
denote $c_1(\sigma)$ by $x$ and $c_1(\tau)$ by $y$, and assume
that $x\geq y$.  Using this terminology, we have $x(x-1)=y(y-1)$.
This equality can hold for integers only if $x,y\in \set{0,1}$, so
$x-y < 2$. However we have $p\mid (x-y)$ for any prime $p\mid m$.
Since $x-y < 2$, we must accept that $x=y$, and so, again by Lemma
\ref{almost-sim}, $\sigma \sim \tau$.
\end{proof}

It's in generalizing to triplets that we first run into trouble.
\begin{theorem}
The representation $\rho_3^\bC$ of $S_n$ on ordered triplets
unites classes if and only if $n$ is even.
\end{theorem}
\begin{proof}
Here $\chi(\sigma) = c_1(\sigma)(c_1(\sigma)-1)(c_1(\sigma)-2)$,
because each fixed point of $T(\sigma)$ comes from an ordered
triplet of fixed points of $\sigma$.  Let $\sigma$ and $\tau$ be a
pair of almost similar elements of order $m$ which get united by
the representation. Again, let $c_1(\sigma)=x$ and $c_1(\tau)=y$
and suppose $x\geq y$. We are facing $x(x-1)(x-2)=y(y-1)(y-2)$.
This can hold for nonnegative integers only if $x, y\in
\set{0,1,2}$, so $x-y < 3$.

If $n$ is even, then consider the following elements:
$\sigma=(1,2)\,(3,4)\,\cdots\,(n-1, n)$ and
$\tau=(1)\,(2)\,(3,4)\,\cdots\,(n-1, n)$.  The orders are prime,
so they are almost similar.  Here $c_1(\sigma)=0$ and
$c_1(\tau)=2$, so $\chi(\sigma)=\chi(\tau)$, and they unite under
$T$.

However, if $n$ is odd, we cannot find such an example.
\begin{itemize}
\item If $m$ has some odd prime component $p$, then as above
$p\mid (x-y)$.  Since $x-y < 3$, we must accept that $x=y$, which
implies (Lemma \ref{almost-sim} again) that $\sigma\sim \tau$.

\item If $m$ is a power of $2$.  In this case both $\sigma$ and
$\tau$ must have at least one fixed point, i.e. $x, y \geq 1$,
because in $n=c_1(\sigma^m) = \sum_{k\mid m, k\neq 1} kc_k(\sigma)
+ c_1(\sigma)$, the left hand side is odd, and all the summands
under the summation sign are even.  Also, we have $2\mid x-y$. So,
if $x-y \neq 0$, then $x-y\geq 2$.  This imples that $x\geq 2+y
\geq 3$, contrary to our assumption about $x$ and $y$. Therefore,
we must agree that $x=y$, and so that $\sigma \sim \tau$.
\end{itemize}
\end{proof}

\begin{theorem}
The representation of $S_n$ on $k$-tuples for $k \geq 4$ always
unites some classes.
\end{theorem}
\begin{proof}
If $n$ is even, then $\sigma=(1,2)\,(3,4)\,\cdots\,(n-1,n)$ and
$\tau=(1)\,(2)\,(3,4)\,\cdots\,(n-1,n)$ furnish an example of
elements that unite, because $c_1(\sigma)=0$ and $c_1(\tau)=2$, so
$\chi(\sigma)=\chi(\tau)$. Similarly, when $n$ is odd:
$\sigma=(1,2)\,(3,4)\,\cdots\,(n-2,n-1)$ and
$\tau=(1)\,(2)\,(3,4)\,\cdots\,(n-2,n-1)$ are united, because
$c_1(\sigma)=1$ and $c_1(\tau)=3$, so $\chi(\sigma)=\chi(\tau)$.
\end{proof}

In summary, we have shown that the natural representation of
$S_n$, i.e., the representation arising from the action of $S_n$
on $1$-tuples, does not unite classes.  The same applies to the
representation arising from the action on pairs.  On triplets it
works only when $n$ is odd, and beyond that, all representations
unite some classes.
\subsection{Representations Arising from the Action of $S_n$ on
Subsets}\label{subsets} The natural route to follow now would be
to try and generalize these results to other permutation
representations, and in particular to those arising from the
action of $S_n$ on $k$-subsets of $\{1,2,...,n\}$. The general
answer eludes us at present, and seems to be pretty
unsatisfactory. However, we have managed to show that the
representation arising from the action of $S_n$ on {\it all}
subsets of $[n]$ does in fact not unite any classes.
\begin{theorem}\label{allsubsets} The action of $S_n$ on the power set $2^{[n]}$ of
$[n]$ does not unite classes.
\end{theorem}
\begin{proof}
As we know, the key is working out the number of fixed points of a
permutation $\sigma$ on this large set. We note that a subset of
$[n]$ is fixed by $\sigma$ if and only if it is a union of cycles
in $\sigma$'s cycle decomposition. This implies that the number of
sets fixed by $\sigma$ is precisely $2^{m(\sigma)}$, where
$m(\sigma)$ is the number of cycles $\sigma$ decomposes into. This
implies that if $\sigma$ and $\tau$ are united by this action,
then $2^{m(\sigma)}=2^{m(\tau)}$, so $m(\sigma)=m(\tau)$. Of
course, the same holds for all powers of $\sigma$ and $\tau$. As
we saw in the proof of Theorem~\ref{similar}, this implies that
$\sigma\sim\tau$.
\end{proof}
We can give a more detailed description: Let us denote by
$\chi_k(\sigma)$ the number of $k$-subsets that $\sigma$ fixes.
Thus $\chi_k(\sigma)$ is the number of ways we can unite some of
$\sigma$'s cycles into $k$-sets. In other words, it's the number
of ways of partitioning $k$, when the available pieces are
$c_1(\sigma)$ ``copies" of 1, $c_2(\sigma)$ ``copies" of 2 etc.
This function can be written out explicitly, but it gets pretty
complicated. For example,
$\chi_2(\sigma)=c_2(\sigma)+{c_1(\sigma)\choose 2}$,
$\chi_3(\sigma)=c_3(\sigma)+c_2(\sigma)c_1(\sigma)+{c_1(\sigma)\choose
3}$. In general,
$$
\chi_k(\sigma)=\sum_{(\lambda_1,...,\lambda_k)\vdash
k}\prod_i{c_i(\sigma)\choose \lambda_i}\,\,\,.
$$
We can write down the generating function
$F_{\sigma}(t)=\sum\chi_k(\sigma)t^k$. It works out to be:
$$
F_{\sigma}(t)=(1+t)^{c_1(\sigma)}(1+t^2)^{c_2(\sigma)}\cdots(1+t^n)^{c_n(\sigma)}\,\,\,.
$$
For example, if $\sigma=(1,2)(3)(4)$, then
$F_{\sigma}(t)=(1+t)^2(1+t^2)=1+2t+2t^2+2t^3+t^4$. This reads that
$\sigma$ fixes only 1 0-subset (hardly surprising...), 2 1-sets
$(\{3\},\{4\})$, 2 2-sets $(\{1,2\},\{3,4\})$, etc...

Consider now the action of $S_n$ on {\it even} sized subsets of
$[n]$. If $n$ is even, then this action unites some classes. For
example, $(1,2)(3,4)...(n-1,n)$ and $(1)(2)(3,4)...(n-1,n)$ get
united. (They are almost similar and both fix $2^{n/2}$ sets.)

However, if $n$ is odd, then this representation does not unite
classes.
\begin{theorem}\label{evensizedsubsets}
Let $n$ be odd. The action of $S_n$ on the set of even-sized
subsets of $[n]$ does not unite classes.
\end{theorem}
\begin{proof}
Let $\sigma$ and $\tau$ be permutations. The cycles of $\sigma$
and $\tau$ cannot all be of even length, because they add up to
$n$, and $n$ id odd. Now plug $-1$ into the generating function
above. On the one hand, the result is 0 because there is an odd
length cycle. On the other hand, it's equal to $\sum_{even
\,\,k}\chi_k(\sigma)-\sum_{odd \,\,k}\chi_k(\sigma)$. Thus
$\sum_{even \,\,k}\chi_k(\sigma)=\sum_{odd
\,\,k}\chi_k(\sigma)=2^{m(\sigma)-1}$. So, if $\sigma$ and $\tau$
are united by the representation, we have $m(\sigma)=m(\tau)$. The
same holds for $\sigma^k$ and $\tau^k$. As we have already
observed, this implies that $\sigma\sim\tau$.
\end{proof}
Finally, we conclude this section by exploring the behavior of the
representation arising from the action of $S_n$ on {\it odd} sized
subsets of $[n]$.
\begin{theorem}\label{oddsizedsubsets}
The action of $S_n$ on the set of odd-sized subsets of $[n]$ does
not unite classes. This does not depend on n's parity.
\end{theorem}
\begin{proof} Using the same reasoning from the previous theorem,
we have that the number of odd-sized subsets of $[n]$ which are
fixed by $\sigma$ is $2^{m(\sigma)-1}$ when $\sigma$ has odd
length cycles. On the other hand, if all $\sigma$'s cycles have
even length, then obviously $\sigma$ cannot fix any odd-sized
subset.

Now, let $\sigma$ and $\tau$ be almost similar permutations, and
suppose that they are united by our representation. This implies
that they fix the same number of odd subsets, and so either both
have some odd cycles, or both have only even sized ones. In the
latter case, they have the same number of 1-cycles, and so by
lemma they are similar. In the former case, we have
$2^{m(\sigma)-1}=2^{m(\tau)-1}$, so $m(\sigma)=m(\tau)$, and as
before, this implies that $\sigma\sim\tau$.
\end{proof}
\subsection{General Fields} The proofs in the previous section
apply only to the complex field $\bC$, (in fact, to all fields
with characteristic $0$.) We shall now show that the same applies
to any field. We shall base ourselves on Theorem \ref{similar}
from Section \ref{sec-1}, where we proved that the natural
representation does not unite classes, regardless the base field.

\begin{lemma}\label{harkava}
Let $f:G\to H$ and $g:H\to K$ be group homomorphisms.
\begin{enumerate}
\item If $f$ and $g$ both do not unite classes, then also $gf$
does not unite them.

\item If $gf$ does not unite classes, then neither does $f$.
\end{enumerate}
\end{lemma}
\begin{proof}\hfill
\begin{enumerate}
\item Suppose that $gf(\sigma)\sim gf(\tau)$.  This implies that
$f(\sigma)\sim f(\tau)$, because $g$ does not unite classes. This
in turn implies that $\sigma\sim \tau$.

\item  Suppose that $f(\sigma)\sim f(\tau)$.  So obviously
$gf(\sigma)\sim gf(\tau)$, and since $gf$ does not unite classes,
we conclude that $\sigma\sim \tau$.
\end{enumerate}
\end{proof}

\begin{theorem}
Let $T$ be any permutation representation of $S_n$.  If $T$ does
not unite classes when considered a representation into $GL(m,
\bC)$, then it does not unite classes when considered as a
representation into $GL(m, \bF)$, for any field $\bF$.
\end{theorem}
\begin{proof}
Any permutation representation can be factored into $S_n\to S_m
\to GL(m, \bC)$, where the first homomorphism is the permutation
representation and the second is the natural representation.  Now,
suppose $T$ does not unite classes.  By Lemma \ref{harkava},
neither does the permutation representation $S_n\to S_m$.  We
already know that the natural representation does not unite
classes, whatever the field. Tacking these two homomorphisms
together gives us the representation in any field, and another
appeal to Lemma \ref{harkava} proves that it still doesn't unite
any classes.
\end{proof}

\section{The action of $S_n \times S_n$ on invertible
matrices}\label{alpha}

In this section we present and application of Theorem
\ref{similar}.
\begin{definition}
Let $\mathbb{F}$ be any field. We define an action of $S_n \times
S_n$ on the group $GL(n,\mathbb{F})$ by
$$
(\pi,\sigma)\bullet A = \pi A
\sigma^{-1}\,\,\,\text{where}\,\,(\pi,\sigma)\in S_n\times S_n
\,\,\text{and}\,\,A\in GL(n,\mathbb F) \eqno(1)
$$
\end{definition}

It is a group action since:
\begin{multline}
(e,e)\bullet A=eAe=A\\
(\pi_1,\sigma_1)\bullet\Big( (\pi_2,\sigma_2)\bullet A\Big)=(\pi_1,\sigma_1)\bullet (\pi_2 A \sigma_{2}^{-1})=\\
\pi_1\pi_2
A\sigma_{2}^{-1}\sigma_{1}^{-1}=(\pi_1\pi_2,\sigma_1\sigma_2)\bullet
A=\Big((\pi_1,\sigma_1)(\pi_2,\sigma_2)\Big)\bullet A\nonumber
\end{multline}

\begin{definition}
Let $M$ be a finite subset of $GL(n,\mathbb{F})$, invariant under
the action of $S_n \times S_n$ defined above. We denote by
$\alpha_M$ the permutation representation of $S_n \times S_n$
obtained from the action $(1)$ . In the sequel we identify the
action $(1)$ with the permutation representation $\alpha_M$
associated with it.
\end{definition}

Now we define a generalization of the conjugacy representation of
$S_n$\label{beta}

 We present a conjugacy representation of $S_n$ on
a subset $M$ of $GL(n,\mathbb{F})$.

\begin{definition}
Denote by $\beta$ the permutation representation of $S_n$ obtained
by the following action on $M$.
$$
\pi\circ A=(\pi,\pi)\bullet A=\pi A\pi^{-1}\eqno(2)
$$
\end{definition}

The connection between $\alpha_M$ and $\beta_M$ is given by the
following easily seen claim:
\begin{claim}
Consider the diagonal embedding of $S_n$ into $S_n\times S_n$.
Then
$$
\beta_M=\alpha_M\downarrow^{S_n\times S_n}_{S_n}.\eqno\qed
$$
\end{claim}

\begin{theorem}\label{2chars}
For every  finite set $M \subseteq GL(n,\mathbb{F})$ invariant
under the action $(1)$ of $S_n \times S_n$ defined above: \item If
$\pi$ and $\sigma$ are conjugate in $S_n$ then
$$ \chi_{\alpha_M}\left(
(\pi,\sigma)\right)=\chi_{\alpha_M}\left(
(\pi,\pi)\right)=\chi_{\beta_M}(\pi)=\#\{A\in M\,|\,\pi
A=A\pi\}\,.
$$
\item If $\pi$ is not conjugate to $\sigma$ in $S_n$ then
$$
\chi_{\alpha_M}\left((\pi,\sigma)\right)=0\,.
$$
\end{theorem}
\begin{proof}
If $\pi$ and $\sigma$ are conjugate in $S_n$ then $(\pi,\sigma)$
is conjugate to $(\pi,\pi)$ in $S_n\times S_n$. Since the
character is a class function, we have:
$$ \chi_{\alpha_M}\big
(\pi,\sigma)=\chi_{\alpha_M} (\pi,\pi)=\#\{A\in M\,|\,\pi
A\pi^{-1}=A\}=\#\{A\in M\,|\,\pi A=A\pi\}
$$
i.e. the value of the character of $\alpha_M$ calculated on the
element $(\pi,\sigma)$ with $\pi$ conjugate to $\sigma$ in $S_n$
is equal to the number of matrices in $M$ which commute with the
permutation matrix $\pi$.

Now, we know that the character of a permutation representation
counts the number of fixed points, so:
$$
\chi_{\alpha_M} (\pi,\sigma)=\#\{A\in M\,|\,\pi
A\sigma^{-1}=A\}=\#\{A\in M\,|\,\pi=A\sigma A^{-1}\}.
$$
Note that $ \pi=A \sigma A^{-1}$ means that $\pi$ and $\sigma$ are
similar as invertible matrices. Thus, by Theorem \ref{similar}, if
$\pi$ and $\sigma$ are not conjugate in $S_n$ they can not be
conjugate in $GL(n,\bF)$ and we have:
$$
\{A\in M\,|\,\pi=A\sigma A^{-1}\}= \varnothing
$$
and so
$$
\chi_{\alpha_M }(\pi,\sigma)=0
$$
if $\pi$ and $\sigma$ are not conjugate in $S_n$.
\end{proof}
This result is applied in~\cite{yona} and~\cite{MW}.

{\bf Acknowledgments.} This paper expands upon the authors theses,
and in particular, the second author's thesis, which was
supervised by Yuval Roichman and Ron Adin. Both authors are very
grateful to them. Also the authors are grateful to Alex Lubotzky,
Eli Bagno, Uzi Vishne and Boris Kunyavskii for helpful
discussions.

\bibliographystyle{amsalpha}

\end{document}